\documentstyle[amsfonts, amssymb, latexsym,12pt]{amsart}

\newtheorem{theorem}{Theorem}[section]
\newtheorem{lemma}[theorem]{Lemma}
\newtheorem{proposition}[theorem]{Proposition}

\newtheorem{definition}[theorem]{Definition}

\def\C{{\mbox{\rm\kern.24em
\vrule width.03em height1.43ex depth-.052ex \kern-.26em C}}}
\def\QSet{\mbox{\rm\kern.24em
\vrule width.03em height1.48ex depth-.051ex \kern-.26em Q}}
\def\Z{{\mbox{\rm\kern.25em
\vrule width.03em height0.57ex depth0ex
\kern.033em
\vrule width.03em height1.52ex depth-0.96ex \kern-.338em Z}}}
\def\Z{{\bf Z}}
\def\N{{{\mbox{\rm I\kern-.2em N}}}}
\def\R{{\mbox{\rm I\kern-.22em R}}}

\def\T{{\bf T}}
\def\supp{{\rm supp}}

\def\dist{{\rm dist}}

\def\111{\gamma}

\def\be#1{\begin{equation}\label{#1}}
\def\bas{\begin{align*}}
\def\eas{\end{align*}}
\def\bi{\begin{itemize}}
\def\ei{\end{itemize}}

\def \endprf{\hfill  {\vrule height6pt width6pt depth0pt}\medskip}
\def\emph#1{{\it #1}}

\title{Multi-parameter paraproducts}

\author{Camil Muscalu}
\address{Department of Mathematics, Cornell University, Ithaca, NY 14853}
\email{camil@@math.cornell.edu}

\author{Jill Pipher}
\address{Department of Mathematics, Brown University, Providence, RI 02912}
\email{jpipher@@math.brown.edu}

\author{Terence Tao}
\address{Department of Mathematics, UCLA, Los Angeles, CA 90095}
\email{tao@@math.ucla.edu}

\author{Christoph Thiele}
\address{Department of Mathematics, UCLA, Los Angeles, CA 90095}
\email{thiele@@math.ucla.edu}

\begin{document}

\begin{abstract}
We prove that the classical Coifman-Meyer theorem holds on any polydisc $\T^d$ of arbitrary dimension $d\geq 1$.
\end{abstract}

\maketitle

\section{Introduction}

This article is a continuation of our previous paper \cite{cjtc}. For $n\geq 1$ let $m (=m(\tau))$ in $L^{\infty}(\R^n)$ be a bounded function,
smooth away from the origin and satisfying

\begin{equation}\label{m1}
|\partial^{\alpha} m(\tau)|\lesssim \frac{1}{|\tau|^{|\alpha|}}
\end{equation}
for sufficiently many multi-indices $\alpha$ \footnote{ $A\lesssim B$ means that there exists an universal constant $C>0$ so that $A\leq CB$.}. Denote by $T_m^{(1)}$ the $n$-linear operator defined by

\begin{equation}\label{tm1}
T_m^{(1)}(f_1,...,f_n)(x) = \int_{\R^n} m(\tau)
\widehat{f_1}(\tau_1)...
\widehat{f_n}(\tau_n) e^{2\pi i x (\tau_1+...+\tau_n)}\, d\tau
\end{equation}
where $f_1,...,f_n$ are Schwartz functions on the real line $\R$. The following statement of Coifman and Meyer is a classical theorem in Analysis \cite{cm},
\cite{ks}, \cite{gt}.

\begin{theorem}\label{cm1}
$T_m^{(1)}$ maps $L^{p_1}\times...\times L^{p_n}\rightarrow L^p$ boundedly, as long as
$1<p_1,...,p_n\leq \infty$, $\frac{1}{p_1}+...+\frac{1}{p_n} = \frac{1}{p}$ and $0<p<\infty$.
\end{theorem}
In \cite{cjtc} we considered the bi-parameter analogue of $T_m^{(1)}$ defined as follows. Let $m( = m(\gamma,\eta))$ in $L^{\infty}(\R^{2n})$ be a bounded function,
smooth away from the subspaces $\{\gamma=0\} \cup \{\eta = 0\}$ and satisfying

\begin{equation}\label{m2}
|\partial^{\alpha}_{\gamma}
\partial^{\beta}_{\eta} 
m(\gamma,\eta)|\lesssim 
\frac{1}{|\gamma|^{|\alpha|}}
\frac{1}{|\eta|^{|\beta|}}
\end{equation}
for sufficiently many multi-indices $\alpha$ and $\beta$.
Denote by $T_m^{(2)}$ the $n$-linear operator defined by

\begin{equation}\label{tm2}
T_m^{(2)}(f_1,...,f_n)(x) = \int_{\R^{2n}} m(\gamma,\eta)
\widehat{f_1}(\gamma_1,\eta_1)...
\widehat{f_n}(\gamma_n,\eta_n) e^{2\pi i x [(\gamma_1,\eta_1)+...+(\gamma_n,\eta_n)]}\, d\gamma d\eta
\end{equation}
where $f_1,...,f_n$ are Schwartz functions on the plane $\R^2$. The following theorem has been proven in \cite{cjtc}.
\begin{theorem}\label{cm2}
$T_m^{(2)}$ maps $L^{p_1}\times...\times L^{p_n}\rightarrow L^p$ boundedly, as long as
$1<p_1,...,p_n\leq \infty$, $\frac{1}{p_1}+...+\frac{1}{p_n} = \frac{1}{p}$ and $0<p<\infty$.
\end{theorem}

The main goal of the present paper is to generalize Theorem \ref{cm2} to the $d$-parameter setting, for any $d\geq 1$.

In general, if $\xi_1= (\xi_1^i)_{i=1}^d,...,\xi_n= (\xi_n^i)_{i=1}^d$ are $n$ generic vectors in $\R^d$, they naturally generate the following
$d$ vectors in $\R^n$ which we will denote by
$\overline{\xi_1}= (\xi_j^1)_{j=1}^n,...,\overline{\xi_d}= (\xi_j^d)_{j=1}^n$. As before, let $m(=m(\xi)=m(\overline{\xi}))$ in
$L^{\infty}(\R^{dn})$ be a bounded symbol, smooth away from the subspaces $\{\overline{\xi_1}=0\}\cup...\cup\{\overline{\xi_d}=0\}$ and satisfying

\begin{equation}\label{m3}
|\partial^{\alpha_1}_{\overline{\xi_1}}...
\partial^{\alpha_d}_{\overline{\xi_d}} m(\overline{\xi})|
\lesssim
\prod_{i=1}^d \frac{1}{|\overline{\xi_i}|^{|\alpha_i|}}
\end{equation}
for sufficiently many multi-indices $\alpha_1,...,\alpha_d$.
Denote by $T_m^{(d)}$ the $n$-linear operator defined by

\begin{equation}\label{tm3}
T_m^{(d)}(f_1,...,f_n)(x) = \int_{\R^{dn}} m(\xi)
\widehat{f_1}(\xi_1)...
\widehat{f_n}(\xi_n) e^{2\pi i x (\xi_1+...+\xi_n)}\, d\xi
\end{equation}
where $f_1,...,f_n$ are Schwartz functions on $\R^d$. The main theorem of the article is the following.

\begin{theorem}\label{cm3}
$T_m^{(d)}$ maps $L^{p_1}\times...\times L^{p_n}\rightarrow L^p$ boundedly, as long as
$1<p_1,...,p_n\leq \infty$, $\frac{1}{p_1}+...+\frac{1}{p_n} = \frac{1}{p}$ and $0<p<\infty$.
\end{theorem}

Classically, \cite{cm}, \cite{ks}, \cite{gt} an estimate as the one in Theorem \ref{cm1} is proved
by using the $T(1)$ theorem of David and Journ\'{e} \cite{stein} together with the Calder\'{o}n-Zygmund decomposition. In particular, the theory
of BMO functions and Carleson measures is involved.

On the other hand, it is well known \cite{cf}, \cite{journe} that in the multi-parameter setting all these results and concepts are much more delicate
(BMO, John-Nirenberg inequality, Calder\'{o}n-Zygmund decomposition). To overcome these difficulties, in \cite{cjtc} we had to develop a completely 
new approach to prove Theorem \ref{cm2}. This approach relied on the one dimensional BMO theory and also on Journ\'{e}'s lemma \cite{journe} \cite{fl}, 
but did not extend to prove the general $d$-parameter case.

The novelty of the present paper is that it simplifies the method introduced in \cite{cjtc} and this simplification works equally well in all dimensions.
Surprisingly, it turned out that one doesn't need to rely on any knowledge of BMO, Carleson measures or Journ\'{e}'s lemma in order to prove the estimates in 
Theorem \ref{cm3}.

We shall rely on our previous paper \cite{cjtc} and for the reader's convenience we chose to present the argument in the same bi-linear bi-parameter setting
(so both $n$ and $d$ will be equal to $2$). However, it will be clear from the proof that its extension to the $n$-linear $d$-parameter case is straightforward.

The paper is organized as follows. In Section 2 we recall the discretization procedure from \cite{cjtc} which reduces the study of our operator to the study of some
general multi-parameter paraproducts. In Section 3 we present the proof of our main theorem, Theorem \ref{cm3} and in the Appendix we give a proof of Lemma \ref{desc} 
which plays 
an important role in our simplified construction.

$\bf{Acknowledgements}$ C.Muscalu and J.Pipher were partially supported by NSF Grants. T.Tao was partially supported by a Packard Foundation Grant.
C.Thiele was partially supported by the NSF Grants 9985572 and DMS 9970469.  Two of the authors (C.Muscalu and C.Thiele) are particularly grateful 
to Centro di Ricerca Matematica Ennio de Giorgi in Pisa for their warm hospitality during their visit in June 2004.

\section{Discrete paraproducts}

As we promised, assume throughout the paper that $n=d=2$. In this case, our operator $T_m^{(d)}$ can be written as

\begin{equation}\label{last}
T_m^{(2)}(f,g)(x)= \int_{\R^4} m(\gamma, \eta)
\widehat{f}(\gamma_1, \eta_1)
\widehat{g}(\gamma_2, \eta_2) e^{2\pi i x[ (\gamma_1, \eta_1)+(\gamma_2, \eta_2)]}\, d\gamma d\eta.
\end{equation}
In \cite{cjtc} Section 1, we decomposed the operator $T_m^{(2)}$ into smaller pieces, well adapted to its bi-parameter structure. This allowed us to reduce its analysis
to the analysis of some simpler discretized dyadic paraproducts. We will recall their definitions below.

An interval $I$ on the real line $\R$ is called dyadic if it is of the form $I=2^k[n,n+1]$ for some $k,n\in\Z$. If $\lambda, t\in [0,1]$ are two parameters
and $I$ is as above, we denote by $I_{\lambda,t}$ the interval $I_{\lambda,t} = 2^{k+\lambda}[n+t, n+t+1]$.

\begin{definition}\label{bump}
For $J\subseteq\R$ an arbitrary interval, we say that a smooth function $\Phi_J$ is a bump adapted to $J$, if and only if the following inequalities hold
\begin{equation}
|\Phi_J^{(l)}(x)|\leq C_{l,\alpha} \frac{1}{|J|^{l}}\frac{1}{\left(1+\frac{\dist(x,J)}{|J|}\right)^{\alpha}},
\end{equation}
for every integer $\alpha\in \N$ and for sufficiently many derivatives $l\in\N$. If $\Phi_J$ is a bump adapted to $J$, we say that $|J|^{-1/2}\Phi_J$ is an
$L^2$ - normalized bump adapted to $J$.
\end{definition}

For $\lambda, t_1, t_2, t_3\in [0,1]$ and $j\in\{1,2,3\}$ we define the discretized dyadic paraproduct $\Pi^j_{\lambda, t_1, t_2, t_3}$  of ``type $j$'' by

\begin{equation}\label{para}
\Pi^j_{\lambda, t_1, t_2, t_3}(f,g)=
\sum_{I\in\cal{D}}\frac{1}{|I|^{1/2}}
\langle f, \Phi^1_{I_{\lambda, t_1}}\rangle
\langle g, \Phi^2_{I_{\lambda, t_2}}\rangle
\Phi^3_{I_{\lambda,t_3}},
\end{equation}
where $f,g$ are complex-valued measurable functions on $\R$ and $\Phi^i_{I_{\lambda, t_i}}$ are $L^2$-normalized bumps adapted to $I_{\lambda, t_i}$
with the additional property that $\int_{\R}\Phi^i_{I_{\lambda, t_i}}(x) dx = 0$ for $i\neq j$, $i=1,2,3$.
$\cal{D}$ is an arbitrary finite set of dyadic intervals and by $\langle\cdot, \cdot\rangle$ we denoted the complex scalar product.

Similarly, for $\vec{\lambda}, \vec{t_1}, \vec{t_2}, \vec{t_3}\in [0,1]^2$ and $\vec{j}\in\{1,2,3\}^2$, we define the discretized dyadic bi-parameter
paraproduct of ``type $\vec{j}$''

$$\Pi^{\vec{j}}_{\vec{\lambda}, \vec{t_1}, \vec{t_2}, \vec{t_3}}= \Pi^{j'}_{\lambda', t'_1, t'_2, t'_3}\otimes
\Pi^{j''}_{\lambda'', t''_1, t''_2, t''_3}$$
by

\begin{equation}\label{para2}
\Pi^{\vec{j}}_{\vec{\lambda}, \vec{t_1}, \vec{t_2}, \vec{t_3}}(f,g) = 
\sum_{R\in\vec{\cal{D}}}\frac{1}{|R|^{1/2}}
\langle f, \Phi^1_{R_{\vec{\lambda}, \vec{t_1}}}\rangle
\langle g, \Phi^2_{R_{\vec{\lambda}, \vec{t_2}}}\rangle
\Phi^3_{R_{\vec{\lambda},\vec{t_3}}},
\end{equation}
where this time $f,g$ are complex-valued measurable functions on $\R^2$, $R=I\times J$ are dyadic rectangles and $\Phi^i_{R_{\vec{\lambda}, \vec{t_i}}}$
are given by

$$\Phi^i_{R_{\vec{\lambda}, \vec{t_i}}} = \Phi^i_{I_{\lambda', t'_i}}\otimes\Phi^i_{J_{\lambda'', t''_i}}$$
for $i=1,2,3$. In particular, if $i\neq j'$ then $\int_{\R}\Phi^i_{I_{\lambda', t'_i}}(x) dx = 0$ and if $i\neq j''$ then
$\int_{\R}\Phi^i_{J_{\lambda'', t''_i}}(x) dx = 0$. $\vec{\cal{D}}$ is an arbitrary finite collection of dyadic rectangles.

We will also denote by $\Lambda^{\vec{j}}_{\vec{\lambda}, \vec{t_1}, \vec{t_2}, \vec{t_3}}(f,g,h)$ the trilinear form given by
\begin{equation}
\Lambda^{\vec{j}}_{\vec{\lambda}, \vec{t_1}, \vec{t_2}, \vec{t_3}}(f,g,h)=
\int_{\R^2}\Pi^{\vec{j}}_{\vec{\lambda}, \vec{t_1}, \vec{t_2}, \vec{t_3}}(f,g)(x,y) h(x, y) dx dy.
\end{equation}
In \cite{cjtc} we showed that Theorem \ref{cm2} can be reduced to the following Proposition.

\begin{proposition}\label{reduction}
Fix $\vec{j}\in \{1,2,3\}^2$ and let $1<p,q<\infty$ be two numbers arbitrarily close to $1$. Let also $f\in L^p$, $\|f\|_p=1$,
$g\in L^q$, $\|g\|_q=1$ and $E\subseteq \R^2$, $|E|=1$. Then, there exists a subset $E'\subseteq E$ with 
$|E'|\sim 1$ such that \footnote{$A\sim B$ means that $A\lesssim B$ and $B\lesssim A$}
\begin{equation}\label{inegalitatea}
\left|\Lambda^{\vec{j}}_{\vec{\lambda}, \vec{t_1}, \vec{t_2}, \vec{t_3}}(f,g,h)\right|\lesssim1
\end{equation}
uniformly in the parameters $\vec{\lambda}, \vec{t_1}, \vec{t_2}, \vec{t_3}\in [0,1]^2$, where $h:=\chi_{E'}$.
\end{proposition}

It is therefore enough to prove the above Proposition \ref{reduction}, in order to complete the proof of our main Theorem \ref{cm2}. Since all the cases are similar, we assume as in
\cite{cjtc} that $\vec{j} = (1,2)$.

To construct the desired set $E'$, we need to recall the ``maximal-square'', ``square-maximal'' and ``square-square'' functions considered in \cite{cjtc}.

For $(x,y)\in\R^2$ define

\begin{equation}\label{ms}
MS(f)(x,y) = \sup_I\frac{1}{|I|^{1/2}}
\left(\sum_{J:R=I\times J\in\vec{\cal{D}}}
\sup_{\vec{\lambda},\vec{t_1}}
\frac{|\langle f, \Phi^1_{R_{\vec{\lambda}, \vec{t_1}}}\rangle|^2}{|J|}
\chi_J(y)
\right)
\chi_I(x),
\end{equation}

\begin{equation}\label{sm}
SM(g)(x,y)= 
\left(
\sum_I
\frac{\sup_{J:R=I\times J\in\vec{\cal{D}}} \sup_{\vec{\lambda},\vec{t_2}}\frac{|\langle g, \Phi^2_{R_{\vec{\lambda}, \vec{t_2}}}\rangle|^2 }{|J|}\chi_J(y) }{|I|}
\chi_I(x)\right)^{1/2}
\end{equation}
and

\begin{equation}\label{ss}
SS(h)(x,y) = 
\left(
\sum_{R\in\vec{\cal{D}}}\sup_{\vec{\lambda}, \vec{t_3}}
\frac{|\langle h, \Phi^3_{R_{\vec{\lambda}, \vec{t_3}}}\rangle|^2 }{|R|}
\chi_R(x,y)
\right)^{1/2}.
\end{equation}
Then, we also recall (see \cite{stein}) the bi-parameter Hardy-Littlewood maximal function

\begin{equation}\label{mm}
MM(F)(x,y) = \sup_{(x,y)\in I\times J}
\frac{1}{|I| |J|}
\int_{I\times J}
|F(x',y')| dx' dy'.
\end{equation}
The following simple estimates explain the appearance of these functions. In particular, we will see that our desired bounds in Theorem \ref{cm2}
can be easily obtained as long as all the indices involved are strictly between $1$ and $\infty$.

We start by recalling the following basic inequality, \cite{cjtc}. If $\Pi^1$ ia a one-parameter paraproduct of ``type $1$'' given by

\begin{equation}
\Pi^1(f_1,f_2) = 
\sum_I
\frac{1}{|I|^{1/2}}
\langle f_1, \Phi^1_I \rangle
\langle f_2, \Phi^2_I \rangle
\Phi^3_I
\end{equation}
then we can write

$$\left|\Lambda^1(f_1,f_2,f_3)\right|=
\left|\int_{\R}\Pi^1(f_1, f_2)(x) f_3(x) dx\right|$$

$$\lesssim
\sum_I
\frac{1}{|I|^{1/2}}
|\langle f_1, \Phi^1_I \rangle|
|\langle f_2, \Phi^2_I \rangle|
|\langle f_3, \Phi^3_I \rangle|$$

$$=
\int_{\R}
\left(
\sum_I
\frac{|\langle f_1, \Phi^1_I \rangle|}{|I|^{1/2}}
\frac{|\langle f_2, \Phi^2_I \rangle|}{|I|^{1/2}}
\frac{|\langle f_3, \Phi^3_I \rangle|}{|I|^{1/2}}
\chi_I(x)
\right)
dx$$

\begin{equation}\label{unu}
\lesssim \int_{\R} M(f_1)(x) S(f_2)(x) S(f_3)(x) dx
\end{equation}
where $M$ denotes the Hardy-Littlewood maximal function and $S$ is the square function of Littlewood and Paley. In particular, we easily see that
$\Pi^1: L^p\times L^q\rightarrow L^r$ for any $1<p,q,r<\infty$ satisfying $1/p+1/q=1/r$.
Analogous estimates hold for any other type of paraproducts $\Pi^j$ for $j=1,2,3$.

Similarly, for the bi-parameter paraproduct $\Pi^{(1,2)}$ of ``type $(1,2)$'' formally defined by $\Pi^{(1,2)} = \Pi^1\otimes \Pi^2$ one obtains the inequalities

$$\left|\Lambda^{(1,2)}(f_1,f_2,f_3)\right|=
\left|\int_{\R}\Pi^{(1,2)}(f_1, f_2)(x,y) f_3(x,y) dx dy\right|$$

\begin{equation}\label{doi}
\lesssim \cdots \lesssim
\int_{\R^2} MS(f_1)(x,y) SM(f_2)(x,y) SS(f_3)(x,y) dx dy,
\end{equation}
and analogous estimates hold for any other type of paraproducts $\Pi^{\vec{j}}$ for $\vec{j}\in \{1,2,3\}^2$.
It is important that all these $MS$, $SM$ and $SS$ functions are bounded on $L^p$ for any $1<p<\infty$. We recall the proof of this fact here (see \cite{cjtc}).
We start with $SM(f_2)(x,y)$. It can be written as

\begin{equation}
SM(f_2)(x,y)= 
\left(
\sum_{\tilde{I}}
\frac{\sup_{\tilde{J}} \frac{|\langle f_2, \Phi^2_{\tilde{I}}\otimes \Phi^2_{\tilde{J}}\rangle|^2 }{|\tilde{J}|}\chi_{\tilde{J}}(y) }{|\tilde{I}|}
\chi_{\tilde{I}}(x)\right)^{1/2}
\end{equation}

$$\lesssim
\left(
\sum_{\tilde{I}}
M(\frac{\langle f_2, \Phi^2_{\tilde{I}}\rangle }{|\tilde{I}|^{1/2}})^2(y)
\chi_{\tilde{I}}(x)\right)^{1/2}
$$
where $\tilde{I}$ and $\tilde{J}$ are the intervals where the corresponding supremums over $\vec{\lambda}, \vec{t_2}\in [0,1]^2$ in (\ref{sm}) are attained.

In particular, by using Fefferman-Stein \cite{fs} and Littlewood-Paley \cite{stein} inequalities, we have

\begin{equation}
\|SM(f_2)\|_p \lesssim
\|
\left(
\sum_{\tilde{I}}
M(\frac{\langle f_2, \Phi^2_{\tilde{I}}\rangle }{|\tilde{I}|^{1/2}})^2(y)
\chi_{\tilde{I}}(x)\right)^{1/2}
\|_p
\end{equation}

$$\lesssim
\|
\left(
\sum_{\tilde{I}}
\frac{|\langle f_2, \Phi^2_{\tilde{I}}\rangle|^2 }{|\tilde{I}|}(y)
\chi_{\tilde{I}}(x)\right)^{1/2}
\|_p\lesssim \|f_2\|_p
$$
for any $1<p<\infty$. Then, we observe that the $MS$ function is pointwise smaller than a certain $SM$ type function and hence bounded on $L^p$, while
the $SS$ function is a classical double square function and its boundedness on $L^p$ spaces is well known, \cite{cf}. As a consequence, it follows as before
that $\Pi^{(1,2)}: L^p\times L^q\rightarrow L^r$ as long as $1<p,q,r<\infty$ with $1/p+1/q=1/r$.

\section{Proof of Proposition \ref{reduction}}
It remains to prove Proposition \ref{reduction}. First, we state the following Lemma.

\begin{lemma}\label{desc}
Let $J\subseteq \R$ be an arbitrary interval. Then, every bump function $\phi_J$ adapted to $J$ can be written as

\begin{equation}
\phi_J = \sum_{k\in\N} 2^{-1000 k} \phi^k_J
\end{equation}
where for each $k\in\N$, $\phi^k_J$ is also a bump adapted to $J$ but with the additional property that $\supp (\phi^k_J)\subseteq 2^k J$
\footnote{$2^k J$ is the interval having the same center as $J$ and whose length is $2^k |J|$.}.
Moreover, if we assume $\int_{R}\phi_J(x) dx = 0$ then all the functions $\phi^k_J$ can be chosen so that $\int_{\R}\phi^k_J(x) dx = 0$
for every $k\in \N$.
\end{lemma}
The proof of this Lemma will be presented later on in the Appendix. It is the main new ingredient which allows us to simplify our previous argument in \cite{cjtc}.
Using it, we can decompose our trilinear form in (\ref{para2}) as

\begin{equation}\label{split}
\Lambda^{\vec{j}}_{\vec{\lambda}, \vec{t_1}, \vec{t_2}, \vec{t_3}}(f,g,h) = \sum_{\vec{k}\in\N^2} 2^{-1000|\vec{k}|}
\sum_{R\in\vec{\cal{D}}}\frac{1}{|R|^{1/2}}
\langle f, \Phi^1_{R_{\vec{\lambda}, \vec{t_1}}}\rangle
\langle g, \Phi^2_{R_{\vec{\lambda}, \vec{t_2}}}\rangle
\langle h, \Phi^{3,\vec{k}}_{R_{\vec{\lambda},\vec{t_3}}}\rangle,
\end{equation}
where the new functions $\Phi^{3,\vec{k}}_{R_{\vec{\lambda},\vec{t_3}}}$ have basically the same structure as the old $\Phi^{3}_{R_{\vec{\lambda},\vec{t_3}}}$
but they also have the additional property that $\supp (\Phi^{3,\vec{k}}_{R_{\vec{\lambda},\vec{t_3}}})\subseteq 2^{\vec{k}}R_{\vec{\lambda},\vec{t_3}}$.
We denoted by $2^{\vec{k}}R_{\vec{\lambda},\vec{t_3}}:= 2^{k_1}I_{\lambda', t'_3}\times 2^{k_2}J_{\lambda'', t''_3}$, $\vec{k}=(k_1, k_3)$ and 
$|\vec{k}|=k_1+k_2$.

Fix now $f,g, E, p, q$ as in Proposition \ref{reduction}. For each $\vec{k}\in\N^2$ define

\begin{equation}
\Omega_{-5|\vec{k}|}=
\{ (x,y)\in\R^2 : MS(f)(x,y)> C2^{5|\vec{k}|} \} \cup \{ (x,y)\in\R^2 : SM(g)(x,y)> C2^{5|\vec{k}|} \}.
\end{equation}
Also, define

\begin{equation}
\tilde{\Omega}_{-5|\vec{k}|}= \{ (x,y)\in\R^2 : MM(\chi_{\Omega_{-5|\vec{k}|}})(x,y) >\frac{1}{100} \}
\end{equation}
and then

\begin{equation}
\tilde{\tilde{\Omega}}_{-5|\vec{k}|}= \{ (x,y)\in\R^2 : MM(\chi_{\tilde{\Omega}_{-5|\vec{k}|}})(x,y) >\frac{1}{2^{|\vec{k}|}} \}.
\end{equation}
Finally, we denote by

$$\Omega = \bigcup_{\vec{k}\in\N^2}\tilde{\tilde{\Omega}}_{-5|\vec{k}|}.$$
It is clear that $|\Omega|< 1/2$ if $C$ is a big enough constant, which we fix from now on. Then, define $E':= E\setminus\Omega$ and observe that
$|E'|\sim 1$. We now want to show that the corresponding expression in (\ref{inegalitatea}) is $O(1)$ uniformly in the parameters
$\vec{\lambda}, \vec{t_1}, \vec{t_2}, \vec{t_3}\in [0,1]^2$. Since our argument will not depend on these parameters, we can assume for simplicity that they are all zero
 and in this case we will write $\Phi^i_R$ instead of $\Phi^i_{R_{\vec{\lambda},\vec{t_i}}}$ for $i=1,2$ and $\Phi^{3,\vec{k}}_R$ instead of
$\Phi^{3,\vec{k}}_{R_{\vec{\lambda},\vec{t_3}}}$.

Fix then $\vec{k}\in\N^2$ and look at the corresponding inner sum in (\ref{split}). We split it into two parts as follows.
Part I sums over those rectangles $R$ with the property that

\begin{equation}
R\cap\tilde{\Omega}_{-5|\vec{k}|}^c \neq \emptyset
\end{equation}
while Part II sums over those rectangles with the property that

\begin{equation}
R\cap\tilde{\Omega}_{-5|\vec{k}|}^c =\emptyset.
\end{equation}

We observe that Part II is identically equal to zero, because if $R\cap\tilde{\Omega}_{-5|\vec{k}|}^c \neq \emptyset$ then $R\subseteq\tilde{\Omega}_{-5|\vec{k}|}$
and in particular this implies that $2^{\vec{k}}R\subseteq \tilde{\tilde{\Omega}}_{-5|\vec{k}|}$ which is a set disjoint from $E'$. It is therefore enough to estimate Part I
only. This can be done by using the technique developed in \cite{cjtc}.

Since $R\cap\tilde{\Omega}_{-5|\vec{k}|}^c\neq\emptyset$, it follows that 
$\frac{|R\cap\Omega_{-5|\vec{k}|}|}{|R|}< \frac{1}{100}$ or equivalently,
$|R\cap\Omega^c_{-5|\vec{k}|}|> \frac{99}{100}|R|$. 

We are now going to describe three decomposition procedures, one for each function $f, g, h$. Later on, we will
combine them, in order to handle our sum.

First, define 

$$\Omega_{-5|\vec{k}|+1}= \{ (x,y)\in\R^2 : MS(f)(x,y)> \frac{C 2^{5|\vec{k}|} }{2^1} \}$$
and set

$$\T_{-5|\vec{k}| +1}= \{ R\in \vec{\cal{D}} : |R\cap\Omega_{-5|\vec{k}| +1}|>\frac{1}{100} |R| \},$$
then define

$$\Omega_{-5|\vec{k}| +2}= \{ (x,y)\in\R^2 : MS(f)(x,y)> \frac{C 2^{5|\vec{k}|} }{2^2} \}$$
and set

$$\T_{-5|\vec{k}| +2}= \{ R\in \vec{\cal{D}}\setminus\T_{-5|\vec{k}| +1} : |R\cap\Omega_{-5|\vec{k}| +2}|>\frac{1}{100} |R| \},$$
and so on. The constant $C>0$ is the one in the definition of the set $E'$ above.
Since there are finitely many rectangles, this algorithm ends after a while, producing the sets $\{\Omega_n\}$
and $\{\T_n\}$ such that $\vec{\cal{D}}=\cup_n\T_n$.

Independently, define

$$\Omega'_{-5|\vec{k}|+1}= \{ (x,y)\in\R^2 : SM(g)(x,y)> \frac{C 2^{5|\vec{k}|} }{2^1} \}$$
and set

$$\T'_{-5|\vec{k}| +1}= \{ R\in \vec{\cal{D}} : |R\cap\Omega'_{-5|\vec{k}| +1}|>\frac{1}{100} |R| \},$$
then define

$$\Omega'_{-5|\vec{k}| +2}= \{ (x,y)\in\R^2 : SM(g)(x,y)> \frac{C 2^{5|\vec{k}|} }{2^2} \}$$
and set

$$\T'_{-5|\vec{k}| +2}= \{ R\in \vec{\cal{D}}\setminus\T'_{-5|\vec{k}| +1} : |R\cap\Omega'_{-5|\vec{k}| +2}|>\frac{1}{100} |R| \},$$
and so on, producing the sets $\{\Omega'_n\}$ and $\{\T'_n\}$ such that $\vec{\cal{D}}=\cup_n\T'_n$.
We would like to have such a decomposition available for the function $h$ also. To do this, we first need to
construct the analogue of the set $\Omega_{-5|\vec{k}|}$, for it. Pick  $N>0$ a big enough integer such that for every
$R\in\vec{\cal{D}}$ we have $|R\cap\Omega^{''c}_{-N}|> \frac{99}{100} |R|$ where we defined

$$\Omega''_{-N}= \{ (x,y)\in\R^2 : SS^{\vec{k}}(h)(x,y)> C 2^N \}.$$
Here $SS^{\vec{k}}$ denotes the same ``square-square'' function defined in (\ref{ss}) but with the functions $\Phi^{3,\vec{k}}_{R_{\vec{\lambda},\vec{t_3}}}$
instead of $\Phi^{3}_{R_{\vec{\lambda},\vec{t_3}}}$
Then, similarly to the previous algorithms, we define

$$\Omega''_{-N+1}= \{ (x,y)\in\R^2 : SS^{\vec{k}}(h)(x)> \frac{C 2^N}{2^1} \}$$
and set

$$\T''_{-N+1}= \{ R\in \vec{\cal{D}} : |R\cap\Omega''_{-N+1}|>\frac{1}{100} |R| \},$$
then define

$$\Omega''_{-N+2}= \{ x\in\R^2 : SS^{\vec{k}}(h)(x)> \frac{C 2^N}{2^2} \}$$
and set

$$\T''_{-N+2}= \{ R\in \vec{\cal{D}}\setminus\T''_{-N+1} : |R\cap\Omega''_{-N+2}|>\frac{1}{100} |R| \},$$
and so on, constructing the sets $\{\Omega''_n\}$ and $\{\T''_n\}$ such that $\vec{\cal{D}}=\cup_n\T''_n$.

Then we write Part I as

\begin{equation}\label{in5}
\sum_{n_1, n_2>-5|\vec{k}|, n_3>-N} \sum_{R\in \T_{n_1, n_2, n_3}}
\frac{1}{|R|^{3/2}}
|\langle f, \Phi^1_{R}\rangle|
|\langle g, \Phi^2_{R}\rangle|
|\langle h, \Phi^{3,\vec{k}}_{R}\rangle| |R|,
\end{equation}
where $\T_{n_1, n_2, n_3}:= \T_{n_1}\cap \T'_{n_2}\cap \T''_{n_3}$. Now, if $R$ belongs to
$\T_{n_1, n_2, n_3}$ this means in particular that $R$ has not been selected at the previous $n_1 -1$, $n_2 -1$ and
$n_3 -1$ steps respectively, which means that $|R\cap\Omega_{n_1-1}|<\frac{1}{100} |R|$,
$|R\cap\Omega'_{n_2-1}|<\frac{1}{100} |R|$ and $|R\cap\Omega''_{n_3-1}|<\frac{1}{100} |R|$ or equivalently, 
$|R\cap\Omega^c_{n_1-1}|>\frac{99}{100} |R|$,
$|R\cap\Omega^{'c}_{n_2-1}|>\frac{99}{100} |R|$ and
$|R\cap\Omega^{''c}_{n_3-1}|>\frac{99}{100} |R|$. But this implies that

\begin{equation}\label{in6}
|R\cap\Omega^c_{n_1-1}\cap\Omega^{'c}_{n_2-1}\cap\Omega^{''c}_{n_3-1}|> \frac{97}{100}|R|.
\end{equation}
In particular, using (\ref{in6}), the term in (\ref{in5}) is smaller than

$$\sum_{n_1, n_2>-5|\vec{k}| , n_3>-N} \sum_{R\in \T_{n_1, n_2, n_3}}
\frac{1}{|R|^{3/2}}
|\langle f, \Phi^1_{R}\rangle|
|\langle g, \Phi^2_{R}\rangle|
|\langle h, \Phi^{3,\vec{k}}_{R}\rangle| |R\cap\Omega^c_{n_1-1}\cap\Omega^{'c}_{n_2-1}\cap\Omega^{''c}_{n_3-1}|=$$

$$\sum_{n_1, n_2> -5|\vec{k}|, n_3>-N} \int_{\Omega^c_{n_1-1}\cap\Omega^{'c}_{n_2-1}\cap\Omega^{''c}_{n_3-1}}
\sum_{R\in \T_{n_1, n_2, n_3}}
\frac{1}{|R|^{3/2}}
|\langle f, \Phi^1_{R}\rangle|
|\langle g, \Phi^2_{R}\rangle|
|\langle h, \Phi^{3,\vec{k}}_{R}\rangle| \chi_{R}(x,y)\, dx dy$$

$$\lesssim \sum_{n_1, n_2> -5|\vec{k}|, n_3>-N} \int_{\Omega^c_{n_1-1}\cap\Omega^{'c}_{n_2-1}\cap\Omega^{''c}_{n_3-1}\cap
\Omega_{\T_{n_1, n_2, n_3}}}
MS(f)(x,y) SM(g)(x,y) SS^{\vec{k}}(h)(x,y)\, dx dy$$

\begin{equation}\label{in7}
\lesssim \sum_{n_1, n_2> -5|\vec{k}|, n_3>-N} 2^{-n_1} 2^{-n_2} 2^{-n_3} |\Omega_{\T_{n_1, n_2, n_3}}|,
\end{equation}
where

$$\Omega_{\T_{n_1, n_2, n_3}}:= \bigcup_{R\in\T_{n_1, n_2, n_3}}R .$$
On the other hand we can write

$$|\Omega_{\T_{n_1, n_2, n_3}}|\leq |\Omega_{\T_{n_1}}|\leq
|\{ (x,y)\in\R^2 : MM(\chi_{\Omega_{n_1}})(x,y)> \frac{1}{100} \}|$$

$$\lesssim |\Omega_{n_1}|= |\{ (x,y)\in\R^2 : MS(f)(x,y)>\frac{C}{2^{n_1}} \}|\lesssim 2^{n_1 p}.$$
Similarly, we have

$$|\Omega_{\T_{n_1, n_2, n_3}}|\lesssim 2^{n_2 q}$$
and also

$$|\Omega_{\T_{n_1, n_2, n_3}}|\lesssim 2^{n_2 \alpha},$$
for every $\alpha >1$. Here we used the fact that all the operators $SM$, $MS$, $SS^{\vec{k}}$, $MM$ are bounded
on $L^s$ (independently of $\vec{k}$) as long as $1<s< \infty$ and also that $|E'|\sim 1$.
In particular, it follows that

\begin{equation}\label{*}
|\Omega_{\T_{n_1, n_2, n_3}}|\lesssim 2^{n_1 p \theta_1}
 2^{n_2 q \theta_2} 2^{n_3 \alpha \theta_3}
\end{equation}
for any $0\leq \theta_1, \theta_2, \theta_3 < 1$, such that $\theta_1+ \theta_2 +\theta_3= 1$.

Now we split the sum in (\ref{in7}) into

\begin{equation}\label{last}
\sum_{n_1, n_2> -5|\vec{k}|, n_3>0} 2^{-n_1} 2^{-n_2} 2^{-n_3} |\Omega_{\T_{n_1, n_2, n_3}}|+ 
\sum_{n_1, n_2> -5|\vec{k}|, 0>n_3>-N} 2^{-n_1} 2^{-n_2} 2^{-n_3} |\Omega_{\T_{n_1, n_2, n_3}}|.
\end{equation}
To estimate the first term in (\ref{last}) we use the inequality (\ref{*}) in the particular case
$\theta_1=\theta_2=1/2$, $\theta_3=0$, while to estimate the second term we use (\ref{*}) for $\theta_j$, $j=1,2,3$
such that $1-p\theta_1>0$, $1-q\theta_2>0$ and $\alpha\theta_3 -1>0$. With these choices, the sum in (\ref{last})
is $O(2^{10|\vec{k}|})$ and this makes the expression in (\ref{split}) to be $O(1)$, after summing over $\vec{k}\in\N^2$.

This completes our proof.

It is now clear that our argument works equally well in all dimensions. In the general case, exactly as in \cite{cjtc} Section 1, one first reduces 
the study of the operator $T_m^{(d)}$ to the study of generic $d$-parameter dyadic paraproducts $\Pi^{\vec{j}}$ for $\vec{j}=(j_1,...,j_d)\in \{1,2,3\}^d$ formally defined by
$\Pi^{\vec{j}} = \Pi^{j_1}\otimes \cdots \otimes \Pi^{j_d}$. Then, one observes as before, by using the linear theory and Fefferman-Stein inequality,
that all the corresponding ``square and maximal'' type functions which naturally appear in
inequalities analogous to (\ref{unu}), (\ref{doi}) are bounded in $L^p$ for $1<p<\infty$ (in fact, as before, it is enough to observe this in the $SS...SMM...M$ case, because
all the other expressions are pointwise smaller quantities).

Having all these ingredients, the argument used in Section 3 works similarly.
 Finally, the $n$-linear case follows in the same way. The details are left to the reader.

\section{Appendix: proof of Lemma \ref{desc}}

In this section we prove Lemma \ref{desc}. Fix $J\subseteq\R$ an interval and let $\phi_J$ be a bump function adapted to $J$. 

Consider $\psi$ a smooth function such that $\supp (\psi) \subseteq [-1/2, 1/2]$ and $\psi = 1$ on $[-1/4, 1/4]$. If $I\subseteq\R$ is a generic interval with center $x_I$,
we denote by $\psi_I$ the function defined by

\begin{equation}
\psi_I(x) = \psi(\frac{x-x_I}{|I|}).
\end{equation}
Since

$$1=\psi_J + (\psi_{2J}-\psi_J) + (\psi_{2^2J}-\psi_{2J})+...$$
it follows that

$$\phi_J = \phi_J\cdot\psi_J + \sum_{k=1}^{\infty}\phi_J\cdot (\psi_{2^kJ}-\psi_{2^{k-1}J})
$$

$$=\phi_J\cdot\psi_J + \sum_{k=1}^{\infty}2^{-1000 k}\cdot[2^{1000 k}\phi_J\cdot (\psi_{2^kJ}-\psi_{2^{k-1}J})]
$$

$$:=\sum_{k=0}^{\infty} 2^{-1000 k} \phi^k_J$$
and it is easy to see that all the $\phi^k_J$ functions are bumps adapted to $J$, having the property that $\supp (\phi^k_J)\subseteq 2^k J$.

Suppose now that in addition we have $\int_{\R}\phi_J(x) dx = 0$. This time, we write

$$\phi_J = \phi_J\cdot\psi_J + \phi_J \cdot(1-\psi_J)$$

$$=\left[\phi_J\cdot\psi_J - \left(\frac{1}{\int_{\R}\psi_J(x) dx}\cdot\int_{\R}\phi_J(x)\psi_J(x) dx\right)\cdot \psi_J\right]$$

$$+ \left[\left(\frac{1}{\int_{\R}\psi_J(x) dx}\cdot\int_{\R}\phi_J(x)\psi_J(x) dx\right)\cdot \psi_J + \phi_J(1-\psi_J)\right]$$

$$:=\phi^0_J + R^0_J.$$
Clearly, by construction we have that $\int_{\R}\phi^0_J(x) dx = 0$ and therefore $\int_{\R}R^0_J(x) dx =0$. Moreover, $\phi^0_J$ is a bump adapted to the interval $J$ 
having the property that $\supp (\phi^0_J)\subseteq J$. On the other hand, since

\begin{equation}
\left|\frac{1}{\int_{\R}\psi_J(x) dx}\cdot\int_{\R}\phi_J(x)\psi_J(x) dx\right|=
\left|\frac{1}{\int_{\R}\psi_J(x) dx}\cdot\int_{\R}\phi_J(x)(1-\psi_J(x)) dx\right|
\end{equation}
$$\lesssim 2^{-1000}$$
it follows that $\|R^0_J\|_{\infty}\lesssim 2^{-1000}$.

Then, we perform a similar decomposition for the ``rest function'' $R^0_J$, but this time we localize it on the larger interval $2J$.
We have

$$R^0_J = R^0_J\cdot\psi_{2J} + R^0_J\cdot(1-\psi_{2J})$$

$$=\left[ R^0_J\cdot\psi_{2J}-\left(\frac{1}{\int_{\R}\psi_{2J}(x) dx}\cdot\int_{\R}R^0_J(x)\psi_{2J}(x) dx\right)\cdot \psi_{2J}\right]$$

$$+ \left[\left(\frac{1}{\int_{\R}\psi_{2J}(x) dx}\cdot\int_{\R}R^0_J(x)\psi_{2J}(x) dx\right)\cdot \psi_{2J} + R^0_J\cdot(1-\psi_{2J})\right]$$

$$:= 2^{-1000}\phi^1_J + R^1_J.$$
As before, we observe that $\int_{\R}\phi^1_J(x) dx = 0$ and also $\int_{\R}R^1_J(x) dx=0$. Moreover, $\phi^1_J$ is a bump adapted to $J$ whose support lies in $2J$
and $\|R^1_J\|_{\infty}\lesssim 2^{-1000\cdot 2}$. Iterating this procedure $N$ times, we obtain the decomposition

\begin{equation}
\phi_J = \sum_{k=0}^N 2^{-1000 k} \phi^k_J + R^N_J
\end{equation}
where all the functions $\phi^k_J$ are bumps adapted to $J$ with $\int_{\R}\phi^k_J(x) dx = 0$ and $\supp (\phi^k_J)\subseteq 2^k J$, while
$\|R^N_J\|_{\infty}\lesssim 2^{-1000 N}$.

This completes the proof of the Lemma.

\end{document}